# DISCUSSION OF "LEAST ANGLE REGRESSION" BY EFRON ET AL.

By Sanford Weisberg[1]

*University of Minnesota*

Most of this article concerns the uses of LARS and the two related methods in the age-old, "somewhat notorious," problem of "[a]utomatic model-building algorithms..." for linear regression. In the following, I will confine my comments to this notorious problem and to the use of LARS and its relatives to solve it.

**1. The implicit assumption.** Suppose the response is $y$, and we collect the $m$ predictors into a vector $x$, the realized data into an $n \times m$ matrix $X$ and the response is the $n$-vector $Y$. If $P$ is the projection onto the column space of $(1, X)$, then LARS, like ordinary least squares (OLS), assumes that, for the purposes of model building, $Y$ can be replaced by $\hat{Y} = PY$ without loss of information. In large samples, this is equivalent to the assumption that the conditional distributions $F(y|x)$ can be written as

$$F(y|x) = F(y|x'\beta) \tag{1.1}$$

for some unknown vector $\beta$. Efron, Hastie, Johnstone and Tibshirani use this assumption in the definition of the LARS algorithm and in estimating residual variance by $\hat{\sigma}^2 = \|(I-P)Y\|^2/(n-m-1)$. For LARS to be reasonable, we need to have some assurance that this particular assumption holds or that it is relatively benign. If this assumption is not benign, then LARS like OLS is unlikely to produce useful results.

A more general alternative to (1.1) is

$$F(y|x) = F(y|x'B), \tag{1.2}$$

where $B$ is an $m \times d$ rank $d$ matrix. The smallest value of $d$ for which (1.2) holds is called the structural dimension of the regression problem [Cook (1998)]. An obvious precursor to fitting linear regression is deciding on the structural dimension, not proceeding as if $d = 1$. For the diabetes data used

[1]Supported by NSF Grant DMS-01-03983.

This is an electronic reprint of the original article published by the Institute of Mathematical Statistics in *The Annals of Statistics*, 2004, Vol. 32, No. 2, 490–494. This reprint differs from the original in pagination and typographic detail.

in the article, the R package dr [Weisberg (2002)] can be used to estimate $d$ using any of several methods, including sliced inverse regression [Li (1991)]. For these data, fitting these methods suggests that (1.1) is appropriate.

Expanding $x$ to include functionally related terms is another way to provide a large enough model that (1.1) holds. Efron, Hastie, Johnstone and Tibshirani illustrate this in the diabetes example in which they expand the 10 predictors to 65 including all quadratics and interactions. This alternative does not include (1.2) as a special case, as it includes a few models of various dimensions, and this seems to be much more complex than (1.2).

Another consequence of assumption (1.1) is the reliance of LARS, and of OLS, on correlations. The correlation measures the degree of linear association between two variables particularly for normally distributed or at least elliptically contoured variables. This requires not only linearity in the conditional distributions of $y$ given subsets of the predictors, but also linearity in the conditional distributions of $a'x$ given $b'x$ for all $a$ and $b$ [see, e.g., Cook and Weisberg (1999a)]. When the variables are not linearly related, bizarre results can follow; see Cook and Weisberg (1999b) for examples. Any method that replaces $Y$ by $PY$ cannot be sensitive to nonlinearity in the conditional distributions.

Methods based on $PY$ alone may be strongly influenced by outliers and high leverage cases. As a simple example of this, consider the formula for $C_p$ given by Efron, Hastie, Johnstone and Tibshirani:

$$(1.3) \qquad C_p(\hat{\mu}) = \frac{\|Y - \hat{\mu}\|^2}{\sigma^2} - n + 2\sum_{i=1}^{n} \frac{\text{cov}(\hat{\mu}_i, y_i)}{\sigma^2}.$$

Estimating $\sigma^2$ by $\hat{\sigma}^2 = \|(I - P)Y\|^2/(n - m - 1)$, and adapting Weisberg (1981), (1.3) can be rewritten as a sum of $n$ terms, the $i$th term given by

$$C_{pi}(\hat{\mu}) = \frac{(\hat{y}_i - \hat{\mu}_i)^2}{\hat{\sigma}^2} + \frac{\text{cov}(\hat{\mu}_i, y_i)}{\hat{\sigma}^2} - \left(\frac{h_i - \text{cov}(\hat{\mu}_i, y_i)}{\hat{\sigma}^2}\right),$$

where $\hat{y}_i$ is the $i$th element of $PY$ and $h_i$ is the $i$th leverage, a diagonal element of $P$. From the simulation reported in the article, a reasonable approximation to the covariance term is $\hat{\sigma}^2 u_i$, where $u_i$ is the $i$th diagonal of the projection matrix on the columns of $(1, X)$ with nonzero coefficients at the current step of the algorithm. We then get

$$C_{pi}(\hat{\mu}) = \frac{(\hat{y}_i - \hat{\mu}_i)^2}{\hat{\sigma}^2} + u_i - (h_i - u_i),$$

which is the same as the formula given in Weisberg (1981) for OLS except that $\hat{\mu}_i$ is computed from LARS rather than from a projection. The point here is that the value of $C_{pi}(\hat{\mu})$ depends on the agreement between $\hat{\mu}_i$ and $\hat{y}_i$, on the leverage in the subset model and on the difference in the leverage between the full and subset models. Neither of these latter two terms has much to do with the problem of interest, which is the study of the conditional distribution of $y$ given $x$, but they are determined by the predictors only.

**2. Selecting variables.** Suppose that we can write $x = (x_a, x_u)$ for some decomposition of $x$ into two pieces, in which $x_a$ represents the "active" predictors and $x_u$ the unimportant or inactive predictors. The variable selection problem is to find the smallest possible $x_a$ so that

$$F(y|x) = F(y|x_a) \tag{2.4}$$

thereby identifying the active predictors. Standard subset selection methods attack this problem by first assuming that (1.1) holds, and then fitting models with different choices for $x_a$, possibly all possible choices or a particular subset of them, and then using some sort of inferential method or criterion to decide if (2.4) holds, or more precisely if

$$F(y|x) = F(y|\gamma' x_a)$$

holds for some $\gamma$. Efron, Hastie, Johnstone and Tibshirani criticize the standard methods as being too greedy: once we put a variable, say, $x^* \in x_a$, then any predictor that is highly correlated with $x^*$ will never be included. LARS, on the other hand, permits highly correlated predictors to be used.

LARS or any other methods based on correlations cannot be much better at finding $x_a$ than are the standard methods. As a simple example of what can go wrong, I modified the diabetes data in the article by adding nine new predictors, created by multiplying each of the original predictors excluding the sex indicator by 2.2, and then rounding to the nearest integer. These rounded predictors are clearly less relevant than are the original predictors, since they are the original predictors with noise added by the rounding. We would hope that none of these would be among the active predictors.

Using the S-PLUS functions kindly provided by Efron, Hastie, Johnstone and Tibshirani, the LARS procedure applied to the original data selects a seven-predictor model, including, in order, BMI, S5, BP, S3, SEX, S6 and S1. LARS applied to the data augmented with the nine inferior predictors selects an eight-predictor model, including, in order, BMI, S5, rBP, rS3, BP, SEX, S6 and S1, where the prefix "r" indicates a rounded variable rather than the variable itself. LARS not only selects two of the inferior rounded variables, but it selects both BP and its rounded version rBP, effectively claiming that the rounding is informative with respect to the response.

Inclusion and exclusion of elements in $x_a$ depends on the marginal distribution of $x$ as much as on the conditional distribution of $y|x$. For example, suppose that the diabetes data were a random sample from a population. The variables S3 and S4 have a large sample correlation, and LARS selects one of them, S3, as an active variable. Suppose a therapy were available that could modify S4 without changing the value of S3, so in the future S3 and S4 would be nearly uncorrelated. Although this would arguably not change the distribution of $y|x$, it would certainly change the marginal distribution of $x$, and this could easily change the set of active predictors selected by LARS or any other method that starts with correlations.

A characteristic that LARS shares with the usual methodology for subset selection is that the results are invariant under rescaling of any individual predictor, but not invariant under reparameterization of functionally related predictors. In the article, the authors create more predictors by first rescaling predictors to have zero mean and common standard deviation, and then adding all possible cross-products and quadratics to the existing predictors. For this expanded definition of the predictors, LARS selects a 15 variable model, including 6 main-effects, 6 two-factor interactions and 3 quadratics. If we add quadratics and interactions first and then rescale, LARS picks an 8 variable model with 2 main-effects, 6 two-factor interactions, and only 3 variables in common with the model selected by scaling first. If we define the quadratics and interactions to be orthogonal to the main-effects, we again get a different result. The lack of invariance with regard to definition of functionally related predictors can be partly solved by considering the functionally related variables simultaneously rather than sequentially. This seems to be self-defeating, at least for the purpose of subset selection.

**3. Summary.** Long-standing problems often gain notoriety because solution of them is of wide interest and at the same time illusive. Automatic model building in linear regression is one such problem. My main point is that neither LARS nor, as near as I can tell, any other *automatic* method has any hope of solving this problem because automatic procedures by their very nature do not consider the context of the problem at hand. I cannot see any solution to this problem that is divorced from context. Most of the ideas in this discussion are not new, but I think they bear repeating when trying to understand LARS methodology in the context of linear regression. Similar comments can be found in Efron (2001) and elsewhere.

SCHOOL OF STATISTICS
UNIVERSITY OF MINNESOTA
1994 BUFORD AVENUE
ST. PAUL, MINNESOTA 55108
USA
E-MAIL: sandy@stat.umn.edu